\documentclass[10pt]{amsart}

\usepackage{amsmath, amssymb, amsthm}

\makeatletter
\def\newrefformat#1#2{%
  \@namedef{pr@#1}##1{#2}}
\def\prettyref#1{\@prettyref#1:}
\def\@prettyref#1:#2:{%
  \expandafter\ifx\csname pr@#1\endcsname\relax%
    \PackageWarning{prettyref}{Reference format #1\space undefined}%
    \ref{#1:#2}%
  \else%
    \csname pr@#1\endcsname{#1:#2}%
  \fi%
}
\makeatother

\newrefformat{eq}{\textup{(\ref{#1})}}
\newrefformat{cha}{Chapter \ref{#1}}
\newrefformat{sec}{Section \ref{#1}}
\newrefformat{tab}{Table \ref{#1} on page \pageref{#1}}
\newrefformat{fig}{Figure \ref{#1} on page \pageref{#1}}

\makeatletter
\def\indsym#1#2{%
  \setbox0=\hbox{$\m@th#1x$}%
  \kern\wd0%
  \hbox to 0pt{\hss$\m@th#1\mid$\hbox to 0pt{$\m@th#1^{#2}$}\hss}%
  \lower.9\ht0\hbox to 0pt{\hss$\m@th#1\smile$\hss}%
  \kern\wd0} \newcommand{\ind}[1][]{\mathop{\mathpalette\indsym{#1}}}
\def\nindsym#1#2{%
  \setbox0=\hbox{$\m@th#1x$}%
  \kern\wd0%
  \hbox to 0pt{\mathchardef\nn="3236\hss$\m@th#1\nn$\kern1.4\wd0\hss}
  \hbox to 0pt{\hss$\m@th#1\mid$\hbox to 0pt{$\m@th#1^{#2}$}\hss}%
  \lower.9\ht0\hbox to 0pt{\hss$\m@th#1\smile$\hss}%
  \kern\wd0}

\makeatother

\def\bb{\bar b} \def\bx{\bar x} \def\by{\bar y} 
 \def\bf{\bar f} \def\bc{\bar c} 
 \def\ba{\bar a}  
\def\bd{\bar d}   
   
    \def\bm{\bar m}

                \def\d{\delta}
         
\def\t{\tau}           \def\ro{\rho}

 \def\P{\mathbb{P}} \def\U{\mathcal{U}}
\def\B{\mathcal{B}} 
   
 \def\A{\mathcal{A}} \def\B{\mathcal{B}}
\def\D{\mathcal{D}}  \def\C{\mathcal{C}} \def\F{\mathcal{F}}
\def\tensor{\otimes} \def\k{\kappa} \def\s{\sigma} 
 \def\L{L}    \def\Sum{\sum}
\def\bra{\langle} \def\ket{\rangle} \def\vphi{\varphi}
  \def\emp{\emptyset} 

\def\F{\mathcal{F}}
\def\I{\mathcal{I}}

\newcommand{\mynewthm}[3][dummythm]{%
  \newtheorem{#2}[#1]{#3}%
  \newrefformat{#2}{#3 \ref{##1}}%
}

\swapnumbers

\theoremstyle{plain}
\mynewthm{theo}{Theorem}
\mynewthm{prop}{Proposition}
\mynewthm{lema}{Lemma}
\mynewthm{fact}{Fact}
\mynewthm{coro}{Corollary}
\mynewthm{ques}{Question}
\mynewthm{exer}{Exercise}

\theoremstyle{definition}
\mynewthm{defi}{Definition}
\mynewthm{rema}{Remark}
\mynewthm{obse}{Observation}
\mynewthm{conv}{Convention}
\mynewthm{nota}{Notation}
\mynewthm{cst}{Construction}
\mynewthm{exam}{Example}

\DeclareMathOperator{\dcl}{dcl}

\DeclareMathOperator{\events}{events}
\DeclareMathOperator{\event}{event}
\DeclareMathOperator{\Step}{Step}

\title{Model theory of probability spaces \\
with an automorphism}

\author{\noindent Alexander Berenstein}
\address{Alexander Berenstein \newline
University of Illinois at Urbana-Champaign \newline
Urbana, Illinois 61801, USA}
\email{aberenst@math.uiuc.edu}

\author{C. Ward Henson}
\address{C. Ward Henson \newline
University of Illinois at Urbana-Champaign \newline
Urbana, Illinois 61801, USA}
\email{henson@math.uiuc.edu}

\thanks{The authors are grateful to Itay Ben-Yaacov for suggesting
that we study the subject of this paper and to Joseph Rosenblatt
and Yevgeny Gordon for helpful conversations.  Research of the
second author was partially supported by NSF grant DMS-0140677}

\begin{document}
\maketitle

\pagestyle{myheadings}
\markboth{\sc \normalsize Alexander Berenstein and C.\ Ward Henson}%
{\sc \normalsize Probability Spaces with an Automorphism}

\begin{abstract}
The class of generic structures among those consisting of the
measure algebra of a probability space equipped with an
automorphism is axiomatizable by positive sentences interpreted
using an approximate semantics.  The separable generic structures
of this kind are exactly the ones isomorphic to the measure
algebra of a standard Lebesgue space equipped with an aperiodic
measure-preserving automorphism. The corresponding theory is
complete and has quantifier elimination; moreover it is stable
with built-in canonical bases.  We give an intrinsic
characterization of its independence relation.
\end{abstract}

\section{Introduction}

One motivation for the work presented in this paper is to
understand generic structures in the setting of measure algebras
equipped with an automorphism. In recent years, existentially
closed (or ``generic'') structures have (again) attracted a lot of
attention in model theory. Several important examples in simple
and stable theories, like algebraically closed fields,
differentially closed fields, random graphs and algebraically
closed fields with a generic automorphism, can be seen as
existentially closed structures. Expanding a stable structure by
adding generic predicates, automorphisms or substructures is a
common way to construct examples of simple theories. It is a
natural idea to test these tools outside the first order context,
in particular to apply them to familiar structures coming from
analysis and probability. The study of some generic expansions of
Hilbert spaces is carried in \cite{BB}. In this paper we examine
measure algebras of probability spaces extended by a generic automorphism.

A second motivation for this paper comes from the more general
objective of understanding the model theory of probability spaces
equipped with an automorphism.  Said differently, it is the
connections between model theory and ergodic theory that are being
explored here, from a very general point of view.  It is not
surprising that some of this general theory might be needed to
study generic structures.  However, it turns out that the
connection between these two projects is much closer than could be
expected.  Indeed, a generic structure in this setting turns out
to be essentially the same as a probability space equipped with an
aperiodic automorphism (\textit{i.e.}, one for which the set of
points with a finite period has measure zero).  Moreover, the
theory of an arbitrary probability space with an automorphism can
be reduced to the theory of its aperiodic part in a simple and
direct way.

The model theoretic background for this paper comes largely from
\cite{HI}, suitably generalized for metric structures.
On a probability space there is a canonical pseudometric, obtained
by taking the distance between sets to be the measure of their
symmetric difference.  The quotient structure which turns this
distance into a metric is exactly the measure algebra of the
original probability space: identify two sets if they differ by a
set of measure zero.  Any automorphism of the probability space
induces one on the measure algebra.  A full discussion of these
structures from a point of view very close to what we take here
can be found in Chapter 7 of \cite{Vlad}.  In general, an
automorphism of a measure algebra need not arise from a point map
of the underlying probability space.  However, if the measure
algebra is separable (as a metric space; equivalently, if it is
countably generated as a measure algebra) then it is isomorphic to
the measure algebra of a Lebesgue space and all of its
automorphisms arise from point maps.

The study of aperiodic maps on Lebesgue spaces is one of the
fundamental parts of ergodic theory. The basics of the theory of
Lebesgue spaces and of aperiodic maps on such spaces is due to
Rokhlin. In section two we introduce the basic definitions and
tools that we need from analysis. Among them is Rokhlin's theorem
on aperiodic maps (stated here as Theorem \ref{rokhlin}), which is
a fundamental tool for this paper.

Stability in probability spaces was studied by Ben-Yaacov in
\cite{BY2}. In section three we present some of the results from
\cite{BY2}, but we translate them to the language of positive
formulas (see \cite{HI}). In particular we show that the
approximate positive theory of atomless probability spaces is well
behaved from the model theoretic point of view: it has quantifier
elimination, is separably categorical and is superstable with
respect to the topology generated by the distance metric in the
space of types. We also prove the definability of this distance
metric and introduce the notion of built-in canonical bases. For
this section we expect the reader to be familiar with the notions
from \cite{HI}, in a modified form suitable to the metric
structures that are studied here.

In section four we construct an axiomatization (denoted by $T_A$)
that isolates the aperiodic expansions among the atomless
probability spaces expanded by an automorphism. We prove that the
class of models of these axioms agree with the existentially
closed models of the approximate theory of probability spaces
expanded by an automorphism. In this section we also show that
$T_A$ has quantifier elimination.

The fifth section is dedicated to showing the stability of
the theory $T_A$ and giving a natural characterization of
non-dividing. We also prove that $T_A$ has built-in canonical
bases.

Entropy is an additive rank defined in ergodic theory. In section six we
review its properties and state some of its model theoretic
consequences.

In section seven we use the properties of entropy to show that the
types in $T_A$ of finite tuples over a set of parameters $B$ are
either non-principal or they belong to the definable closure of
$B$. Using an argument like one that has been applied to ACFA, as
well as the properties of entropy, we generalize a theorem
of ergodic theory and show that the types of transformally
independent elements are orthogonal to the types of transformally
definable elements (see Definition \ref{transformally}).

Finally in the last section we discuss the model theory of an
arbitrary (not necessarily atomless) probability space extended
by a general (not necessarily aperiodic) automorphism. We show
the approximate positive theory of such a structure is easily
reduced to the theory of its aperiodic part.

\section{Probability spaces, Lebesgue spaces and aperiodic maps}
In this section we present the basic information about probability
spaces including Lebesgue spaces, their measure algebras of
events, and aperiodic maps.

A \emph{probability space} is a triple $(X,\B,P)$, where $X$ is a
space, $\B$ is a $\s$-algebra of subsets of $X$ and $P$ is a
measure on $\B$ such that $P(X)=1$. We say that a probability
space $(X,\B,P)$ is \emph{atomless} if for any $B\in \B$ such that
$P(B)>0$, there are $B_1,B_2\in \B$ such that $B=B_1 \cup B_2$,
$B_1$ and $B_2$ are disjoint and $P(B_1)>0$, $P(B_2)>0$.

We say that $A_1,A_2\in \B$ determine the same \emph{event}, if the
symmetric difference of the sets, denoted by$A_1\triangle A_2$,
has measure zero. We denote the class of $A$ under this
equivalence relation by $a$. Throughout this section lower case
letters will stand for events and capital letters either for
elements of the $\s$-algebra or for sets of events.

The operations of complement, union and intersection are well
defined for events. We call $\events(X,\B,P)$ a \emph{measure algebra} and the pair $(\events(X,\B,P),P)$ a \emph{measured algebra}.

Whenever $\C\subset \B$ is a $\s$-subalgebra, we denote by
$\events(\C,P)$ the measure algebra of events coming from $\C$
with respect to $P$. Whenever $P$ is clear from context, we just
write $\events(\C)$. Conversely, whenever $C\subset
\events(X,\B,P)$, we denote by $\bra C\ket$ the
$\sigma$-subalgebra of $\B$ generated by the elements $\{A\in \B:
a\in C\}$.

There are two approaches to understand isomorphisms on probability
spaces. On the one hand, we have point maps between the spaces, on
the other, measure preserving maps between measured algebras.

\begin{defi}
Let $(X_1,\B_1,m_1)$, $(X_2,B_2,m_2)$ be probability spaces and
let $\hat \B_1=\events(X_1,\B_1,m_1)$, $\hat
\B_2=\events(X_2,\B_2,m_2)$ be their measure algebras. By an
\emph{isomorphism} of the measured algebras we mean a bijection
$\Phi\colon\hat \B_2\to \hat B_1$ which preserves complements,
countable unions and intersections and satisfies $m_1(\Phi(\hat
B))=m_2(\hat B)$ for all $\hat B\in \hat \B_2$. The probability
spaces are said to be \emph{conjugate} if their measured algebras
are isomorphic.
\end{defi}

\begin{defi}
Let $(X_1,\B_1,m_1)$, $(X_2,B_2,m_2)$ be probability spaces and
let $\hat \B_1=\events(X_1,\B_1,m_1)$, $\hat
\B_2=\events(X_2,\B_2,m_2)$ be their measure algebras. Let 
$M_1\in \B_1$, $M_2\in \B_2$ with $m_1(M_1)=1=m_2(M_2)$. An
invertible measure preserving transformation $\phi\colon M_1\to
M_2$ is called an \emph{isomorphism} between $(X_1,\B_1,m_1)$ and
$(X_2,B_2,m_2)$. If $(X_1,\B_1,m_1)=(X_2,B_2,m_2)$, we call $\phi$
an \emph{automorphism}. The induced map $\phi\colon \hat B_1\to 
\hat B_2$ is called an \emph{induced isomorphism} of the measured 
algebras.
\end{defi}

To bridge the gap between these two approaches we need to know how
point maps are related to maps of measured algebras. Clearly any
two isomorphic probability spaces are conjugate; however, the
converse does not hold in general.  The next definition concerns
a well-known special class of probability spaces which are well
behaved from this point of view.

\begin{defi}
A probability space $(X,\B,m)$ is a \emph{Lebesgue space} if it is
isomorphic to a probability space which is the disjoint union of a
countable (or finite) set of points $\{y_1,y_2,\dots\}$, each of
positive measure, and the space $([0,s],\L([0,s]),l)$, where
$\L([0,s])$ is the Lebesgue $\sigma$-algebra of $[0,s]$ and $l$ is
Lebesgue measure. Here $s=1-\Sum_{i=1}^{\infty} p_i$, where
$p_i>0$ is the measure of $\{y_i\}$.
\end{defi}

The theory of Lebesgue spaces was developed by Rokhlin. On these
spaces the notion of isomorphism and conjugacy coincide (see
Theorem 2.2 in \cite{Wal}). Thus, as long as we work on Lebesgue
spaces, we can switch between point maps and maps on the measured
algebra of events.

For the rest of this section we fix $(X,\B,m)$ an atomless Lebesgue 
space. It is shown in \cite{Ha} that for any $A,B\in \B$ such that 
$m(A)=m(B)$, there is an automorphism $\eta$ of the space such that
$\mu(\eta(A)\triangle B)=0$.

Let $G$ be the group of measure preserving automorphisms on
$(X,\B,m)$, where we identify two maps if they agree on a
set of measure one. There is a natural representation of $G$ in
$\B(L^2(X,\B,m))$ (the space of bounded linear operators on
$L^2(X,\B,\mu)$); it sends $\tau\in G$ to the unitary operator
$U_\tau$ defined for all $f\in L^2(X,\B,\mu)$ by $U_\tau(f)=f
\circ \tau$. The norm topology on $\B(L^2(X,\B,\mu))$ pulls back
to a group topology on $G$, which is called in \cite{Ha} the 
\emph{uniform topology} on $G$. For $\tau,\eta\in G$, let
$\ro(\tau,\eta)=m(\{x\in X: \tau(x)\neq \eta(x)\})$. It is shown
in \cite{Ha} that $\ro$ is a metric for the uniform topology.

For the rest of this section we will study aperiodic maps and
their properties. A good source for this material is the book of
Halmos \cite{Ha} on ergodic theory.

\begin{defi}
Let $(Y,\C,\mu)$ be an atomless probability space and let $\tau$
be an automorphism of $(Y,\C,\mu)$. We say that $\tau$ is
\emph{aperiodic} if for every $n\in \mathbb{N}^+$, the set $\{x\in
X: \tau^{n}(x)=x\}$ has measure zero.
\end{defi}

One of the key tools in studying aperiodic automorphisms is the
following theorem by Rokhlin:

\begin{theo}(Rokhlin's Lemma \cite{Ha},\cite{Sh})
\label{rokhlin}
Let $(Y,\C,\mu)$ be an atomless probability space and $\tau$ an
aperiodic automorphism of this space. Then for every positive
integer $n$ and $\epsilon>0$, there exist a measurable set $E$
such that the sets $E,\tau(E),\dots,\tau^{n-1}(E)$ are disjoint and
$\mu(\cup_{i<n} \tau^i(E))> 1-\epsilon$.
\end{theo}

\begin{defi}
We call a map $\eta \in G$ a \emph{cycle} of period $k$ if there
is a set $E\in \B$ such that $E$,\dots,$\eta^{k-1}(E)$ forms a
partition of $X$ and $\eta^k=id$.
\end{defi}

\begin{obse}\label{unifapprox}
Let $\tau\in G$ be aperiodic. By Rokhlin's Lemma, for every $N>0$
there is a cycle $\eta \in G$ of period $N$ such that
$\ro(\tau,\eta)\leq 2/N$ (see \cite[pp. 75]{Ha}).
\end{obse}

\begin{rema}\label{conjugacy} Any two cycles $\eta_1,\eta_2\in G$ of period $k$ are
conjugate in $G$. Let $E\in \B$ be such that
$E$,\dots,$\eta_1^{k-1}(E)$ forms a partition of $X$ and let $F\in
\B$ be such that $F$,\dots,$\eta_2^{k-1}(F)$ forms a partition of
$X$. Since $(X,\B,\mu)$ is a Lebesgue space, there is a measure
preserving invertible map $\gamma$ such that $\gamma(E)=F$.
Extend $\gamma$ by defining for $x=\eta_1^i(y)\in \eta^i(E)$, 
$\gamma(x)=\eta_2^i(\gamma(y))$. Then $\gamma \eta_1=\eta_2 \gamma$.
\end{rema}

\begin{prop}\label{isom}
Let $\tau_1,\tau_2\in G$ be aperiodic. Then for every
$\epsilon>0$, there is a conjugate $\tau_2'$ of $\tau_2$ such that
$\ro(\tau_1,\tau_2')\leq \epsilon$.
\end{prop}

\begin{proof}
Let $N>0$ be such that $4/N<\epsilon$. By \ref{unifapprox}, we
can find $\eta_1,\eta_2\in G$ cycles of period $N$ such that
$\ro(\tau_i,\eta_i)< 2/N$ for $i=1,2$. Using \ref{conjugacy} we get
$\gamma\in G$ such that $\eta_1=\gamma^{-1}\eta_2\gamma$. Let 
$\tau_2'=\gamma^{-1}\tau_2 \gamma$. Then $\ro(\tau_1,\tau_2')\leq
\ro(\tau_1,\eta_1)+\ro(\eta_1,\tau_2')=\ro(\tau_1,\eta_1)+\ro(\eta_2,\tau_2)<4/N$.
\end{proof}

\section{The model theory of probability spaces}

We develop the model theory of probability spaces inside
structures of the form
$M=(\events(X,\B,P),\emp,X,^{-1},\cap,\cup,P)$, where $(X,\B,P)$
is an atomless probability space, $\emp$ is the event
corresponding to $\emp$, $X$ is the event corresponding to $X$,
$\cup$, $\cap$ stand for the union and intersection of events,
$^{-1}$ for the complement of events and $P$ for the probability
of events. We also include in $M$ a second sort for the ordered
field of real numbers and constants for all rationals. For $a,b\in
M$, let $\ro(a,b)=P(a \triangle b)$. The distance $\ro$ is a
metric on the space of events, it is definable from $P$ and makes
$M$ a complete metric space. We will use the tools from \cite{HI}, modified
for metric structures, to understand the model theory of probability
spaces (and later their expansion by generic automorphisms). We
call a structure $M$ as above a \emph{probability structure}.

We define \emph{positive formulas} inductively. If $q\in
\mathbb{Q}$ and $x_1,\dots,x_n$ are variables in the sort of events
and $t(y_1,\dots,y_n)$ is a polynomial with coefficients in
$\mathbb{Q}$, then $t(P(x_1),\dots,P(x_n))\geq q$ and
$t(P(x_1),\dots,P(x_n))\leq q$ are positive formulas. If $\vphi$,
$\psi$ are positive formulas, so are $\vphi\vee \psi$ and $\vphi
\wedge \psi$. Finally, if $\vphi$ is a formula, so are $\exists x
\vphi$ and $\forall x \vphi$, where $x$ is a variable in the sort
of events.

From an abstract point of view, the structures considered here
consist of a complete metric space $(M,\rho)$ equipped with
operations making $M$ a Boolean algebra on which $P(a) =
\rho(a,0)$ defines $P$ to be a probability measure and it is 
translation invariant under the operation of symmetric 
difference. In Chapter 7 of \cite{Vlad} there is a full discussion 
of the fact that these structures are exactly the measured algebras 
of probability spaces.

Strictly speaking \cite{HI} is formulated in the setting of normed
space structures, and the probability structures considered here
are not of that type.  However, the aspects of \cite{HI} on which
we rely here are routinely seen to apply to probability
structures, and we will cite results from \cite{HI} (such as the
existence of highly saturated and homogeneous models) without
additional comment.  Note that since a probability structure is
based on a bounded metric space (of diameter 1) there is no need
to bound quantifiers; the key aspects of \cite{HI} that are
essential to what we do here are the use of positive formulas
(only) and the use of an approximate semantics.

In this section we denote by lower case letters the events and by
capital letters elements in the $\s$-algebra and sets of events.

We need the following special case of the Radon-Nikodym theorem:

\begin{theo}(Radon-Nikodym)
Let $(X,\B,P)$ be an atomless probability space, let $\C\subset
\B$ be a $\s$-subalgebra and let $A\in \B$. Let $a$ be the event
corresponding to $A$. Then there is a unique $g_a\in L^1(X,\C,P)$
such that for any $B\in \C$, $\int_{B}g_a dP=\int_{B} \chi_{A}
dP$. Such an element $g_a$ is called the \emph{conditional
probability of $a$ with respect to $\C$} and it is denoted by
$\P(a|\C)$.
\end{theo}

\begin{defi}
Let $\kappa$ be a regular cardinal larger that $2^{\aleph_0}$ . We
say that a metric structure $M$ in a language $L$ is a
\emph{$\kappa$-universal domain} if it is $\kappa$-strongly
homogeneous and $\kappa$-saturated for all reducts of the language
(see \cite{HI}). We call a subset $C\subset M$ \emph{small} if
$|C|<\kappa$.
\end{defi}

\begin{defi} Let $M$ be a metric structure which is a
$\kappa$-universal domain for its positive theory. Let $C\subset
M$ be small and $n\in \mathbb{N}^+$. By an \emph{$n$-type}
over $C$ we mean the collection of positive formulas with
parameters in $C$ realized by some tuple $\bar a\in M^n$.
\end{defi}

The collection of $n$-types over a small set $C$ is independent of
the choice of $M$; it only depends on $Th_{\A}(M,c)_{c\in C}$.

We denote by $T$ the approximate theory \cite{HI} of an atomless
probability structure.

\begin{defi} Let $M$ be a metric structure which is a
$\kappa$-universal domain for its theory. Let $C\subset M$ be
small. The \emph{definable closure of $C$}, denoted by $\dcl(C)$,
is the collection of elements in $M$ that are fixed under
automorphisms of $M$ fixing $C$ pointwise.
\end{defi}

The following two lemmas are proved in \cite{BY2}:

\begin{lema}\label{types}
Let $M\models T$ be a $\kappa$-universal domain, let $C\subset M$
be small and let $a_1,\dots,a_n;b_1,\dots,b_n\in M$. Then
$tp(a_1,\dots,a_n/C)=tp(b_1,\dots,b_n/C)$ iff \newline
$\P(a_1^{i_1}\wedge\dots\wedge a_n^{i_n} | \bra C \ket
)=\P(b_1^{i_1}\wedge\dots\wedge b_n^{i_n} |\bra C \ket)$ for all
$i_j\in \{1,-1\}$, $j\leq n$.
\end{lema}

In particular, any probability structure has quantifier
elimination (see \cite[pp.~86--88]{HI}).

\begin{lema}\label{dcl}
Let $M\models T$ be a $\kappa$-universal domain and let $C\subset
M$ be small. Then $\dcl(C)=\events(\bra C\ket)$.
\end{lema}

The next lemma gives another basic fact about the model theory of probability spaces.

\begin{lema}
The positive theory $T$ is separably categorical.
\end{lema}

\begin{proof}
Let $N$ be a separable complete model of $T$. Then $N$ is a probability
structure coming from an atomless Lebesgue space $(X,\B,m)$, where
$\B$ is countably generated. Hence $N$ is isomorphic to the
probability structure coming from the standard interval
$([0,1],\B,m)$.
\end{proof}

Let $M$ be a metric structure which is a $\kappa$-universal
domain for its theory and denote by $\ro$ the metric of $M$. Then
the collection of $n$-types over small sets also form a metric
space \cite{HI}. When $C\subset M$ is small and $p,q$ are $n$-types over
$C$, we define the \emph{$d$-metric} by $d(p,q)=\inf\{\max_{i\leq
n} \ro(a_i,b_i): (a_1,\dots,a_n)\models p, (b_1,\dots,b_n)\models
q\}$.

In general, in a metric space structure the distance between
types over finite sets need not be definable. The next lemma shows
that these expressions on probability structures are actually
definable (a fact known to analysts). We sketch a proof; see also
lemma 6.3 in \cite{Sh}:

\begin{lema}\label{distance}
Let $M\models T$ be a $\kappa$-universal domain and let $C\subset
M$ be small. Let $\ba=(a_1,\dots,a_n)\in M^n$, $\bb=(b_1,\dots,b_n)\in
M^n$ be partitions of the probability structure. Then
$d(tp(\ba/C),tp(\bb/C))=\max_{i\leq n}\|\P(a_i|\bra C\ket
)-\P(b_i|\bra C\ket)\|_1$, where $\| \ \ \|_1$ is the $L_1$-norm.
\end{lema}

\begin{proof}
We prove the result for $C=\emp$. It is sufficient to show that
$d(tp(\ba),tp(\bb))\leq Max_{i\leq n}|P(a_i)-P(b_i)|$, since the other
direction is obvious. Since $T$ is separably categorical and
$C=\emp$, we can work in the probability structure of
the standard Lebesgue space $([0,1],\B,P)$. Let $(A_1,\dots,A_n)\in
\B^n$, $(B_1,\dots,B_n)\in \B^n$ be partitions of $[0,1]$.
Reordering the partitions if necessary, we may assume there is
$m\leq n$ such that $P(A_i)\geq P(B_i)$ iff $i\leq m$. Let $A_i'$
be the closed interval $[\sum_{j<i}P(A_j),\sum_{j\leq i}P(A_j)]$
for $i\leq n$. Write $a_i,b_i,a_i'$ for the events determinded by 
$A_i$, $B_i$ and $A_i'$ for $i\leq n$. Then
$tp(a_1,\dots,a_n)=tp(a_1',\dots,a_n')$. For $i\leq m$, let $B_i'$ be
the closed interval $[\sum_{j<i}P(A_j),\sum_{j<i}P(A_j)+P(B_i)]$.
Clearly $P(A_i'\triangle B_i')=P(A_i)-P(B_i)$ for $i\leq m$. For
$i>m$, let $B_i'$ be the set $[\sum_{j<i}P(A_j),\sum_{j\leq
i}P(A_j)]\cup D_i'$, where $\{D_i':m<i\leq n\}$ are measurable
sets, disjoint from each other and disjoint from
$B_1',\dots,B_m',A_{m+1}',\dots,A_n'$, such that
$P(B_i)=P(A_i)+P(D_i')$ for $m<i\leq n$. Let $b_i'$ be the event determined by $B_i'$.
Then $tp(b_1,\dots,b_n)=tp(b_1',\dots,b_n')$ and $P(a_i'\triangle b_i')=
|P(a_i)-P(b_i)|$ for $i\leq n$.
\end{proof}

Note that the previous lemma strengthens Lemmas \ref{types} and \ref{dcl}.

\begin{defi}
Let $M$ be a metric structure which is a $\kappa$-universal
domain for its theory.. Let $C\subset M$ be small and let $(\I,<)$
be a countable infinite linear order. Let $I=(\ba_i:i\in \I)$ be a
sequence of $n$-tuples. We say that $I$ is \emph{indiscernible
over $C$} if for any $m\in \mathbb{N}$ and elements
$i_1<i_2<\dots<i_{2m}$ of $\I$,
$tp(\ba_{i_1},\dots,\ba_{i_m}/C)=tp(\ba_{i_{m+1}},\dots,\ba_{i_{2m}}/C)$.
\end{defi}

Let $M\models T$ be a $\kappa$-universal domain, let $C\subset M$
be small and let $I=(a_i:i\in \omega)$ be an indiscernible
sequence of $1$-tuples over $C$. Let $I=(\chi_{a_i}:i\in \omega)$
be the corresponding sequence of characteristic functions.
Then the sequence $I'$ is spreadable \cite[pp.168]{Ka}.
Furthermore, for all $k_1<\dots<k_n\in \omega$ and $i_1,\dots,i_n \in
\{-1,1\}$, we have $\P(a_1^{i_1}\wedge\dots.\wedge a_n^{i_n}|\bra
C\ket)=\P(a_{k_1}^{i_1}\wedge\dots.\wedge a_{k_n}^{i_n}|\bra
C\ket)$.

\begin{defi}
Let $M$ be a metric structure which is a $\kappa$-universal
domain for its theory. Let $\bd_0\in M^m$ and let $\ba \in M^n$.
We say that $tp(\ba/C\cup \bd_0)$ \emph{does not divide over $C$}
if for any indiscernible sequence $I=(\bd_i:i<\omega)$ over $C$,
there is $\ba'\in M^n$ such that
$tp(\ba',\bd_i/C)=tp(\ba,\bd_0/C)$ for all $i<\omega$. Let
$D\subset M$ be small. We say that $tp(\ba/C\cup D)$ does not
divide over $C$ if for all finite $\bd\subset D$, $tp(\ba/C\cup
\bd)$ does not divide over $C$. Whenever $tp(\ba/C\cup D)$ does
not divide over $C$ we say that $\ba$ is \emph{independent from
$D$ over $C$} and we write $\ba \ind_{C} D$. We say that $tp(\ba/C)$ 
is \emph{stationary} if whenever $\bb\in M^n$ and $D\supset C$ is small,
 $tp(\ba/C)=tp(\bb/C)$ and $\ba,\bb \ind_{C} D$ implies that $tp(\ba/D)=tp(\bb/D)$.

Let $\ba_0\in M^n$ and let $I=(\ba_i:i\in \omega)$ be an
indiscernible sequence over $C$. We say that $I$ is a \emph{Morley
sequence} in $tp(\ba_0/C)$ if for every $n\in \mathbb{N}$, 
$tp(\ba_{n+1}/C\cup \ba_0\dots\cup \ba_{n})$ does not divide 
over $C$. 
\end{defi}

We refer the reader to \cite{Io1,Io2,BY0,BY3} for the properties
of non-dividing (non-forking) and stable structures.

In \cite[Theorem 2.10]{BY2} there is a natural characterization of
non-dividing in probability structures:

\begin{prop}
Let $M\models T$ be a $\kappa$-universal domain and let $C\subset
M$ be small. Let $a_1,\dots,a_n;b_1,\dots,b_m\in M$. Then
$tp(a_1,\dots,a_n/C\cup \{b_1,\dots,b_m\})$ does not divide over $C$
if and only if $\P(a_1^{i_1}\wedge\dots\wedge a_n^{i_n}|\bra C\cup
\{b_1,\dots,b_m\}\ket )=\P(a_1^{i_1}\wedge\dots\wedge a_n^{i_n}|\bra C
\ket)$ for all $i_j\in \{-1,1\}$, $j\leq n$.
\end{prop}

\begin{prop}
Let $M\models T$ be a $\kappa$-universal domain. Then
\begin{enumerate}
\item $M$ is stable and types over sets are stationary.
\item $M$ is $\omega$-stable with respect to the $d$-metric.
\end{enumerate}
\end{prop}

\begin{proof}
(1) It is proved in \cite[section 2.3]{BY2}.
(2) Let $C\subset M$ be countable. We may assume that $C$
is closed under finite intersections, unions and complements. Let
$\Step(C)$ be the set of step functions in $L^1(X,\bra C\ket,m)$ with
coefficients in $\mathbb{Q}$ and let $\F=\{tp(a/C): \P(a|\bra C
\ket)\in \Step(C)\}$. Then $\F$ is a countable set of types. By
Lemma \ref{distance}, $\F$ is a dense subset of the space of
$1$-types over $C$ with respect to the $d$-metric. Then by
\cite{Io1}, $M$ (equivalently $T$) is $\omega$-stable with respect
to the $d$-metric.
\end{proof}

\begin{rema}\label{superstable}
In particular, $M$ is superstable with respect to the $d$-metric:
for any $\ba\in M^n$, $C\subset M$ small and $\epsilon>0$, there
is $C_0\subset C$ finite and $\ba'\in M^n$ such that $P(a_i
\triangle a_i')<\epsilon$ for $i=1,\dots,n$ and ${\ba'}\ind_{C_0}C$.
The proof follows along the same lines as the proof of
part (2) in the previous proposition.

Let $\epsilon>0$ and let $B\subset C\subset M$ be small sets. We
say that $tp(a/C)$ \emph{$\epsilon$-divides over $B$} if
$d(tp(a/C),tp(a'/BC))\geq \epsilon$, where $tp(a'/C)$ is the
(unique) non-dividing extension of $tp(a/B)$. Let $SU_{\epsilon}(tp(a/B))$ be
the foundation rank of $\epsilon$-dividing of the type $tp(a/B)$. Then
for any $\epsilon>0$, $a\in M$ and $B\subset M$ small,
$SU_{\epsilon}(tp(a/B))$ is finite. This property translates the 
condition of being superstable of finite $SU$-rank into the current 
metric setting.
\end{rema}

We now recall some definitions from \cite{BB}. These definitions
apply to general metric structures (not just probability
structures).

\begin{defi}
Let $M$ be a metric structure which is a $\kappa$-universal
domain for its theory. Let $I,J\subset (M^n)^\omega$ be countable
indiscernible sequences. We say $I$ and $J$ are \emph{colinear} if
the concatenation $IJ$ is an indiscernible sequence. We say that
$I$ and $J$ are \emph{parallel} if there is another infinite
indiscernible sequence $K\subset (M^n)^\omega$ such that $I$ is
colinear to $K$ and $J$ is colinear to $K$.
\end{defi}

Let $M$ be a metric space structure which is a $\kappa$-universal
domain for its theory and assume that $Th_{\A}(M)$ is stable. We
will show that parallelism is an equivalence relation. Assume that
$I_1$ is parallel to $I_2$ and that $I_2$ is parallel to $I_3$. So
there are $K_1=(a_i:i<\omega)$ and $K_2=(b_i:i<\omega)$
indiscernible sequences such that all concatenations $I_1K_1$,
$I_2K_1$, $I_2K_2$, $I_3K_2$ are indiscernible. Let
$tp(a/I_2I_1I_3K_1K_2)$ be the unique non-dividing extension of
$tp(a_1/I_2)$. Then by stability, $tp(a/I_2I_1I_3K_1K_2)$ is also
the non-dividing extension of $tp(b_1/I_2)$. Let $L$ be a Morley
sequence in $tp(a/I_1I_2I_3K_1K_2)$. Then $I_1L$ and $I_3L$ are
colinear, so $I_1$ is parallel to $I_3$.

\begin{defi}\label{defCb} Let $M$ be a metric structure
which is a $\kappa$-universal domain for its theory. Let $C\subset
M$ be small, let $\ba_0\in M^n$ and assume that $tp(\ba_0/A)$ is
stationary. We say that $tp(\ba_0/A)$ has a \emph{built-in
canonical base} if there is a small set $B\subset M$ such that for
some (equivalently, any) Morley sequence $I=(\ba_i: i<\omega)$ over $A$, 
the parallelism class of $I$ is interdefinable with $B$. That is, for
every automorphism $\Psi$ of $M$, $B$ is fixed (pointwise) by $\Psi$
 if and only if the parallelism class of $I$ is fixed (setwise) 
by $\Psi$. We call $B$ a \emph{built-in canonical base for} 
$tp(\ba_0/A)$. We say that $M$ has \emph{built-in canonical bases} if 
for all $C\subset M$ small, $n\in \mathbb{N}$ and $\ba_0\in M^n$ such 
that $tp(\ba_0/A)$ is stationary, there is a built-in canonical base
for $tp(\ba_0/A)$.
\end{defi}

\begin{rema}\label{probCb} Let $M\models T$ be a $\kappa$-universal
domain, let $C\subset M$ be small and let $a_1,\dots,a_n\in M$. Let
$\D$ be the (smallest) measure algebra making \newline
$\P(a_1^{i_1}\wedge\dots\wedge a_n^{i_n}|\bra C\ket )$ measurable
for all $i_j\in \{-1,1\}$, $j\leq n$. Then $\events(\D)$ is a
built-in canonical base for $tp(a_1,\dots,a_n/C)$. An exposition of
this fact from a slightly different perspective can be found in
\cite{BY2} or \cite{BBH}.
\end{rema}

\section{Aperiodic algebras and quantifier elimination}

We start by fixing the notation for the rest of the paper. Recall
that we denote by $L$ be the language of probability structures
and by $T$ the approximate theory of atomless probability
structures. Write $L_{\tau}$ for the language $L$ expanded by a
unary function with symbol $\tau$ and let $T_{\tau}$ be the theory
$T\cup ``\tau$ is an automorphism''.

Let $(X,\B,\mu)$ be an atomless Lebesgue space and let $M$ be the
probability structure of $(X,\B,\mu)$. Let $G$ be the
group of automorphisms of this space, where we identify two maps
if they agree on a set of measure one. Let $\tau, \ro\in G$. Note
that the map sending $A\in \B$ to $\ro^{-1}(A)$ is a measure
preserving automorphism that induces an isomorphism between the
structures $(M,\tau)$, $(M,\ro^{-1}\tau\ro)$. Let $\tau_1,
\tau_2\in G$ be aperiodic. An application of proposition
\ref{isom} with the values for $\epsilon$ ranging over the
sequence $\{1/n:n\in \mathbb{N}^+\}$, shows that there are
countable ultrapowers of $(M,\tau_1)$ and $(M,\tau_2)$ which are
isomorphic. Thus any two aperiodic transformations in a
\emph{Lebesgue space} have the same approximate elementary theory.
The aim of this section is to study the approximate theory of
probability structures expanded by an aperiodic point
automorphism. Since the elements of a probability structure
are events, we need to define a notion of aperiodicity with
respect to the measure algebra of events.

\begin{defi}
Let $(X,\B,m)$ be a probability space, let $\hat \B$ be the
corresponding measure algebra of events and let $\tau$ be an
automorphism of the measure algebra $\hat B$. The map $\tau$ is
called \emph{aperiodic} if for all
$n\in \mathbb{N^+}$ and $\epsilon>0$ there is $b\in \hat \B$ 
such that $m(b\cap\tau^n(b))\leq \epsilon$ and $|m(b)-1/2|\leq \epsilon$.
\end{defi}

Note that the previous definition can be expressed by positive
formulas in the language $L_{\tau}$. Denote by $T_A$ the theory
$T_{\tau}\cup ``\tau$ is aperiodic''. $T_A$ is $T_\tau$ plus a set 
of approximations, for every $n\in \mathbb{N}$, of the existence 
of a set $b$ of measure $1/2$ which is disjoint from $T^n(b)$,

\begin{lema}\label{equiv}
Let $(X,\B,m)$ be a probability space and let $\hat \B$ be the
corresponding measure algebra of events. Let $\tau$ be an
automorphism of $(X,\B,m)$. Then $\tau$ is aperiodic iff the 
induced automorphism $\tau$ on the measured algebra 
$(\hat \B,m)$ is aperiodic.
\end{lema}

\begin{proof}
If $\tau$ is an aperiodic automorphism of $(X,\B,m)$, then Rokhlin's 
Lemma the induced automorphism $\tau$ on the measured algebra of events 
is also aperiodic.

Now assume that the induced automorphism on the measured algebra of 
events is aperiodic. By a way of contradiction, assume that $\tau$ is not
aperiodic automorphism of $(X,\B,\mu)$. Then for some $n>0$, 
$m(\{x:\tau^n(x)=x\})\geq \epsilon$ for some $\epsilon>0$. Since $\tau$ 
is an automorphism, $\{x:\tau^n(x)=x\}\in \B$. Let $a$ be the event 
corresponding to the set $\{x:\tau^n(x)=x\}$ and let $\delta>0$ be such 
that $9 \delta<\epsilon$. Since $\tau$ is aperiodic with respect to the
measure algebra of events, there is an event $b$ such that
$m(b\cap \tau^n(b))\leq \delta$ and $|m(b)-1/2|\leq \delta$.
Then $m(a)=m(a\cap b)+m(a\cap \tau^n(b)\cap b^c)+m(a\cap b^c\cap
\tau^n(b^c))$. Observe that $m(b^c\cap \tau^n(b^c))=
1-m(b)-m(\tau^n(b))+m(b\cap \tau^n(b))\leq 1-2(1/2-\delta)+\delta=3\delta$.
Since $m(a)>9\delta$, we must have that $m(a\cap b)>3\delta$ or that
$m(a\cap \tau^n(b))>3\delta$.
Assume that $m(a\cap b)>3\delta$. Then $a\cap b\subset a$ and thus
$\tau^n(a\cap b)=a\cap b$. But $m(b\cap \tau^n(b))\leq \delta$ which
is a contradiction. A similar contradiction can be found when
$m(a\cap \tau^n(b))>3\delta$. Thus $\tau$ has to be aperiodic.

\end{proof}

The previous Lemma shows that the two notions of aperiodicity coincide 
for probability structures expanded by point automorphisms and that
the approximate positive theory of atomless Lebesgue space
expanded by an aperiodic automorphism is axiomatized by $T_A$.

\begin{lema}
$T_A$ is complete.
\end{lema}

\begin{proof}
Let $(M_1,\tau_1)$ and $(M_2,\tau_2)$ be two models of $T_A$. Then
there are separable complete models $(M_i',\tau_i')\models T_A$ which are
elementarily equivalent to $(M_i,\tau_i)$ for $i=1,2$
respectively. By separable categoricity of $T$, for $i=1,2$,
$M_i'$ is isomorphic to the probability structure associated to a
 Lebesgue space and $\tau_i'$ can be assumed to
be a point automorphism of the underlying probability space. By
the previous lemma, $\tau_1'$ and $\tau_2'$ are aperiodic
automorphisms as point set maps and thus by Proposition \ref{isom}
we have $(M_1',\tau_1')\equiv_{\A} (M_2',\tau_2')$.
\end{proof}

\begin{defi} Let $M$ be a metric structure. We say that $Th_{\A}(M)$
is \emph{model complete} if for all $N_0\subset N_1$ models of
$Th_{\A}(M)$, $N_0\preceq_{\A} N_1$.
\end{defi}

Our next aim is to show that $T_A$ is model complete. The
techniques used for the proof are similar to the ones used in
proving the completeness of the theory $T_{A}$, but now we need to
include parameters.

\begin{lema}\label{conjugate}
Let $(M,S)$ be a separable complete model of $T_{\tau}$, where $S$ is an
$(n+1)$-shift. Let $f_1,\dots,f_m\in M$. Let $(N,\s)\supset
(M,S)$ be a separable complete extension which is a model of $T_{\tau}$ and
in which $\s$ is an $(n+1)$-shift. Then there is $\colon:M\to N$
measure preserving such that $\eta S=\s \eta$ and $\eta(f_k)=f_k$
for $k\leq m$.
 \end{lema}

\begin{proof}
We may assume that $\{f_1,\dots,f_m\}$ are the atoms of an algebra
and that $S(\{f_1,\dots,f_m\})=\{f_1,\dots,f_m\}$. Since $N$ is
separable, we may assume that there is a Lebesgue space
$(X,\B,\mu)$ such that $N=\events(\B)$ and that $\s$ is induced by a point
map (also denoted by $\s$) such that $\s^{n+1}=id$. We may also suppose
that $M=\events(\C)$ for some complete subalgebra $\C$ of $\B$ and 
that $S$ is induced by a point map $S$ such that $S^{n+1}=id$.

Let $A\in \C$ be such that $(A,\dots,S^n(A))$ is a partition (up to 
measure zero) of $X$. Let $F_1,\dots,F_m\in \C$ be disjoint sets such that 
$\event(F_i)=f_i$ for $i\leq m$. Consider $\eta$ a measure
preserving bijection between $(A,\C\upharpoonright_{A},m)$ and
$(A,\B\upharpoonright_{A},m)$ such that $\eta(A\cap F_k)=A\cap
F_k$ for $k\leq m$. Since $F_1,\dots,F_m$ are the atoms of an
algebra, $\eta$ is well defined for all $x\in A$. Extend $\eta$ by
defining, for $x\in A\cap F_k$, $\eta(S^{i}(x))=\s^{i}(\eta(x))$
for $i\leq n$. Then $\eta(S^{i}(A\cap F_k))=\s^{i}(A\cap F_k)$ for
$i\leq n$, $k\leq m$ and thus $\eta(f_k)=f_k$ for $k\leq m$. By 
the definition of $\eta$ we also have $\eta S=\s \eta$.
\end{proof}

\begin{prop}
$T_A$ is model complete.
\end{prop}

\begin{proof}
Let $(M,\tau)\models T_{A}$ be separable complete and let 
$(N,\tau_1)\models T_{\tau}$ be a separable complete extension of 
$(M,\tau)$. Note that $(N,\tau_1)\models T_{A}$. Let 
$\bf=(f_1,\dots,f_m)\in M^m$.

By compactness we can find a separable complete elementary
extension $(M_1,\tau)$ of $(M,\tau)$ such that for every $n>0$ there is
$a\in M_1$ such that $(a,\dots,\tau^n(a))$  forms a partition of
$M_1$. We can amalgamate $(M_1,\tau)$ and $(N_1,\tau_1)$ over
$(M,\tau)$ (see the discussion of \emph{relative independent joinings over a common factor} in \cite[Chapter 6]{Gl}). Call this structure $(N_1,\tau_1)$. To
prove the proposition it suffices to show that any formula
$\vphi(\bx,\bf)$ true in $(N_1,\tau_1)$ and any approximation
$\vphi'(\bx,\bf)$ of $\vphi(\bx,\bf)$, there is a realization of
$\vphi'(\bx,\bf)$ in $(M_1,\tau)$.

Since $(N_1,\tau_1)$ is separable, we may assume that there is an
atomless Lebesgue space $(X,\B,\mu)$ such that
$N_1=\events(X,\B,\mu)$ and $\tau_1$ is an automorphism of
$(X,\B,\mu)$. Furthermore, we may assume that there is a complete
$\s$-subalgebra $\C$ of $\B$ such that $M_1=\events(\C)$ and that 
$\tau$ is an isomorphism of $(X,\C,\mu)$.

For every $n>0$, we will construct a measure preserving isomorphism 
$\eta_n\colon (X,\C,\mu)\to (X,\B,\mu)$ such that $\eta_n(f_k)=f_k$
for $k\leq m$ and $\ro_{M_1}(\eta_n^{-1} \tau_1 \eta_n, \tau)\leq
2/(n+1)$.

We start by finding approximations of $\tau$ and $\tau_1$ by 
$(n+1)$-shifts. Let $A\in \C$ with event $a$. For 
$x\in \cup_{i<n}\tau_1^i(A)$, let $\s_1(x)=\tau_1(x)$ and for
$x\in \tau_1^n(A)$, let $\s_1(x)=\tau_1^{-n}(x)$. Then $\s_1^{n+1}=id$, 
$(a,\dots,\s_1^n(a))$ is a partition of $M_1$ and 
$\ro(\tau_1,\s_1)\leq 1/(n+1)$. For 
$x\in \cup_{i<n}\tau^i(A)$, let $\s(x)=\tau(x)$ and for
$x\in \tau^n(A)$, let $\s(x)=\tau^{-n}(x)$. Then $\s^{n+1}=id$ and
$\ro(\tau,\s)\leq 1/(n+1)$. 

By Lemma \ref{conjugate} there is a measure preserving isomorphism 
$\eta_n\colon(X,\C,\mu)\to (X,\B,\mu)$ such that 
$\eta_n \s=\s_1 \eta_n$ and
$\eta_n(f_k)=f_k$ for $k\leq m$. Note that $\ro(\eta_n^{-1} \tau_1
\eta_n,\tau)\leq 2/(n+1)$.

Since $n$ was arbitrary, we can find countable ultrapowers
$(M_2,\tau)$ and $(N_2,\tau)$ of $(M_1,\tau)$ and $(N_1,\tau)$
respectively such that $(M_2,\tau,\bf)\cong (N_2,\tau,\bf)$ and 
hence $(M_1,\tau,\bf)\equiv_{\A}(N_1,\tau,\bf)$.
\end{proof}

\begin{defi}
We say that $(M,\tau)\models T_\tau$ is \emph{existentially
closed} if whenever $(N,\t)\supset (M,\t)$, $\bm\in M^n$ and 
$\vphi(\bx,\bm)$ is a quantifier free formula such that 
$(N,\t)\models \exists \bx \vphi(\bx,\bm)$, then for any
approximation $\vphi'(\bx,\by)$ of $\vphi(\bx,\by)$,
$(M,\t)\models \exists \bx \vphi'(\bx,\bm)$.
\end{defi}

\begin{lema}
The models of $T_A$ are precisely the existentially closed models
of $T_{\tau}$.
\end{lema}

\begin{proof}
Since any model of $T_\tau$ can be extended to a model of $T_A$,
an existentially closed model of $T_\tau$ is a model of $T_A$. The
other direction follows from the previous proposition.
\end{proof}

\begin{rema}
The authors initially studied this subject to answer a question of
Itay Ben-Yaacov about the axiomatizability of probability spaces
expanded by generic automorphisms. Indeed, $T_A$ is an
axiomatization for this class. Another axiomatization comes from a
suggestion of Anand Pillay. It is given by the following axioms
indexed by $\epsilon\in \mathbb{Q}^+$ and the arities of $\bx$ and $\ba$:

\begin{enumerate}
\item[(ECN)] $\forall \bx \forall \by \forall \ba \ \big{(} d( tp(\tau(\bx)/\tau(\ba)),tp(\by/\tau(\ba)))\geq \epsilon  \vee \\
d(tp(\bx'/\ba\tau(\ba)),tp(\bx/\ba \tau(\ba)))\geq \epsilon  \vee \\
d(tp(\by'/\ba\tau(\ba)),tp(\by/\ba \tau(\ba)))\geq \epsilon \vee \\
\exists \bc \
d(tp(\bc,\tau(\bc)/\ba,\tau(\ba)),tp(\bx,\by/\ba,\tau(\ba)))\leq
2 \epsilon$.
\end{enumerate}

where $d$ is the distance between types and
$tp(\bx'/\bd,\tau(\bd))$, $tp(\by'/\bd,\tau(\bd))$ are the unique
non-dividing extensions of $tp(\bx/\bd)$ and $tp(\by/\tau(\bd))$
respectively.

This axiomatization says that any possible extension of an
automorphism $\tau$ is approximately realized already in $\tau$.
It is an exercise to the reader to show that this scheme
axiomatizes the existentially closed models of $T_\tau$.

The advantage of this approach is that it can also be used to show
the existence of a model companion for other structures, for
example Hilbert spaces expanded by an automorphism.
\end{rema}

\begin{obse}
The theory ACFA in characteristic $0$ can be seen as the limit,
as the characteristic $p$ goes to infinity, of the theories of
algebraically closed fields expanded by adding the Frobenius
automorphism. The theory of aperiodic automorphisms on atomless
Lebesgue spaces is the limit, as $n$ goes to infinity, of the
theory of a probability space formed by $n$ points of equal weight
with a cycle of period $n$.
\end{obse}

Let $(M,\tau)$ be a $\kappa$-universal domain of $T_A$. We denote
the definable closure in the structure $M$ by $\dcl$ and the
definable closure in the structure $(M,\tau)$ by $\dcl_{\tau}$. It
is easy to give, as in first order theories (see \cite{CP}), a
characterization for $\dcl_\tau$:

\begin{lema} Let $(M,\tau)$ be a $\kappa$-universal domain of $T_A$ and let
$\ba\subset M$. Then $\dcl_{\tau}(\ba)=\dcl(\{\tau^i(\ba): i\in
\mathbb{Z}\})$.
\end{lema}

In ergodic theory, \emph{joinings} give different ways of
amalgamating two probability structures with automorphisms into a
common extension. In particular, the \emph{relative independent
joining over a common factor} (described in Section 6.1 of
\cite{Gl}) corresponds to the model-theoretic \emph{free amalgamation}.
Since $T_A$ is model complete and has the amalgamation property,
it should have quantifier elimination. That is the content of the
next theorem:

\begin{theo}
$T_A$ has elimination of quantifiers.
\end{theo}

\begin{proof}
Let $(M,\tau)\models T_A$ and let $\ba,\bb\in M^n$ such that
$qftp_\tau(\ba)=qftp_\tau(\bb)$. Then for every $k<\omega$,
$tp(\tau^{-k}(\ba),\dots,\tau^{k}(\ba))=tp(\tau^{-k}(\bb),\dots,\tau^k(\bb))$.
Let $f$ be an $L$-isomorphism taking $(\tau^i(\ba):i\in
\mathbb{Z})$ to $(\tau^i(\bb):i\in \mathbb{Z})$. The map $f$ has a
unique extension from $\dcl_\tau(\ba)=\dcl(\tau^i(\ba):i\in
\mathbb{Z})$ to $\dcl_\tau(\bb)$, which is an $L$-isomorphism. By
stationarity of types in atomless probability spaces, the sets
$\dcl_\tau(\ba)$ and $\dcl_\tau(\bb)$ are amalgamation bases in
$T_\tau$. Since $(M,\tau)$ is existentially closed by a back and
forth argument $f$ is an $L_\tau$-isomorphism and
$tp_\tau(\ba)=tp_\tau(\bb)$.
\end{proof}

\section{Independence and stability}

In this section we introduce an abstract notion of independence
and show that it agrees with non-dividing. This idea follows the
approach used in \cite{CP} to characterize non-dividing inside a
first order stable structure expanded by a generic automorphism.
We reserve the use of the word independence for independence of
events in the sense of probability structures. Fix
$(M,\tau)\models T_A$ a $\kappa$-universal domain.

\begin{defi}
Let $\ba\in M^n$ and let $C\subset B\subset M$ be small. We
say that $\ba$ is $\tau$-independent from $B$ over $C$ and
write $\ba \ind[\tau]_{C} \B$ if $\dcl_{\tau}(\ba)$ is
independent from $\dcl_{\tau}(B)$ over $\dcl_{\tau}(C)$.
\end{defi}

The next lemma shows that types in $T_A$ are stationary with
respect to $\tau$-independence. The main tool for this proof is
quantifier elimination for $T_A$.

\begin{prop}
Let $\ba,\bb\in M^n$ and let $C\subset D\subset M$.
Suppose that $tp_{\tau}(\ba/C)=tp_{\tau}(\bb/C)$ 
and that $\ba \ind[\tau]_{C} D$ and $\bb \ind[\tau]_{C} \D$. Then 
$tp_{\tau}(\ba/D)=tp_{\tau}(\bb/D)$.
\end{prop}

\begin{proof}
Let $\ba,\bb,C,D$ be as above. Then for every $k<\omega$,
$$tp(\tau^{-k}(\ba),\dots,\tau^k(\ba))/\dcl_{\tau}(C))=tp(\tau^{-k}(\bb),\dots,\tau^k(\bb))/\dcl_{\tau}(C)).$$
By stationarity of types in probability spaces, we get
$tp(\tau^{-k}(\ba),\dots,\tau^k(\ba))/\dcl_{\tau}(D))=tp(\tau^{-k}(\bb),\dots,\tau^k(\bb))/\dcl_{\tau}(D))$.
Since this equality holds for all $k<\omega$, by quantifier
elimination of $T_A$, $tp_\tau(\ba/D)=tp_\tau(\bb/D)$.
\end{proof}

\begin{coro}
The theory $T_A$ is stable and $\tau$-independence agrees with
non-dividing.
\end{coro}

\begin{proof}
By the properties of independence in $M$, it is clear that
$\tau$-independence satisfies: symmetry, transitivity, extension,
local character and finite character (see \cite{Wa} for the
definition of these properties). By the previous proposition it
also satisfies stationarity. Since non-dividing can be
characterized by these properties, $\tau$-independence agrees with
non-dividing (see \cite{Wa}, \cite{BY3}) and $T_A$ is stable.
\end{proof}

\begin{rema}
Since $T$ has built-in canonical bases (see Definition
\ref{defCb}, Remark \ref{probCb}), $T_A$ will also have built-in
canonical bases.  For any $c_1,\dots,c_m\in M$ and $A\subset M$,
denote by $Cb(c_1,\dots,c_m/A)$ a built-in canonical base for
$tp(c_1,\dots,c_m/A)$. Let $\ba=(a_1,a_2,\dots,a_n)\in M^n$ and let
$A\subset M$ be such that $A=\dcl_\tau(A)$. Then a built-in
canonical base for $tp_\tau(\ba/A)$ is $\cup \{Cb(c_1,\dots,c_m/A):
c_1,\dots,c_m\in \dcl_{\tau}(a_1,\dots,a_n), m<\omega\}$. This is
reminiscent of what happens in ACFA, see \cite{Ma,Ch}.
\end{rema}

Roughly speaking, stability as developed in \cite{Io1,Io2,Io3}
corresponds to the study of universal domains that have a bound on
the size of the space of types. This analysis is carried out
through the density character of uniform structures. Independence
is studied through the notion of non-forking and stability turns
out to be equivalent to definability of types (see section 3 in
\cite{Io2}). The analysis of independence developed in \cite{BY3}
is based on the notion of non-dividing (defined by Shelah). A
structure is stable when it has definability of types (see section
2 in \cite{BY3}) and inside a stable structure, non-dividing can
be characterized by its properties. Hence both points of view
coincide and furthermore, non-forking in \cite{Io1,Io2,Io3}
corresponds to non-dividing from \cite{BY3}.

The previous corollary shows that $T_A$ is stable. Now we will
explicitly count types. We will show that $T_A$ is $\omega$-stable
with respect the \emph{minimal uniform structure}, which was
introduced in \cite{Io1}. We recall the definition:

\begin{defi}
Let $(M,\tau)\models T_A$ be a $\kappa$-universal domain and let
$B\subset M$ be small. Given $L_\tau$ formulas
$\phi_1(\bx,\by)<\phi_1'(\bx,\by)$,\dots,$\phi_k(\bx,\by)<\phi_k'(\bx,\by)$,
define $U[\phi_1,\phi_1',\dots.,\phi_k,\phi_k']$ as the set of all
pairs of $\tau$-types $(p,q)$ with parameters in $B$ such that
$\phi_i(\bx,\bb)\in p$ implies $\phi'_i(\bx,\bb)\in q$ and
$\phi_i(\bx,\bb)\in q$ implies $\phi'_i(\bx,\bb)\in p$, for
$i=1,\dots,k$ and $\bb \subset B$.

The family $U[\phi_1,\phi_1',\dots.,\phi_k,\phi_k']$ forms a uniform
structure on $S_n(A)$ and it is called the \emph{minimal uniform
structure}.
\end{defi}

Before we show that $T_A$ is $\omega$-stable with respect to the
minimal uniform structure, we need to introduce some new
definitions.

\begin{defi}
Let $(M,\tau)\models T_A$, let $A\subset M$ and let $b\in M$. We
say that $b$ is \emph{$m$-step independent over $A$} if
$tp(b/A\cup \{\tau^j(b):j\geq 1\})$ does not divide over $A\cup
\{\tau^j(b):1\leq j\leq m\}$.
\end{defi}

\begin{defi}
Let $(M,\tau)\models T_A$, let $A\subset M$ be a subalgebra and
let $b\in M$. We say that $b$ is \emph{$m$-step simple over $A$}
if for any $c\in \dcl(\tau^{i}(b): 0\leq i\leq m)$, $\P(c|\bra
A\ket)$ is a finite sum of rational multiples of characteristic
functions of elements of $A$.
\end{defi}

\begin{theo}
$T_A$ is $\omega$-stable with respect to the minimal uniform
structure.
\end{theo}

\begin{proof}
Let $(M,\tau)\models T_A$ be a $\k$-universal domain and let
$A_0\subset M$ be countable. We need to prove that there is a
countable dense subset of the space of $\tau$-types over $A_0$
with respect to the minimal uniform structure. Let $A_1=\cup_{i\in
\mathbb{Z}} \tau^{i}(A_0)$, so $A_1$ is countable. Finally let $A$
be the boolean algebra generated by $A_1$. Note that $A$ is
countable. To show $\omega$-stability, it is enough to find a
countable dense subset of the space of $\tau$-types over $A$ with
respect to the minimal uniform structure. Let $\F$ be the set of
all $\tau$-types $tp_\tau(b/A)$ where there is $m$ such that $b$
is m-step simple over $A$ and $m$-step independent over $A$. The
set $\F$ is countable.

\textbf{Claim:} $\F$ is dense in the space of types over $A$ with
respect to the minimal uniform structure.

Let $\vphi(x,\by)<\vphi'(x,\by)$ be $L_\tau$-formulas. Then, by
quantifier elimination, there are $m<\omega$ and
$\psi(x_1,\dots,x_m,\by_1,\dots,\by_m)<\psi'(x_1,\dots,x_m,\by_1,\dots,\by_m)$
quantifier free $\L$-formulas such that $(M,\tau)\models
\vphi(x,\by)\implies \psi(x,\dots,\tau^m(x),\by,\dots,\tau^m(\by))$
and $(M,\tau)\models
\psi'(x,\dots,\tau^m(x),\by,\dots,\tau^m(\by))\implies \vphi'(x,\by)$.

By the perturbation lemma, there is $\epsilon>0$ such that
whenever $$(\bc,\bd)\models
\psi(\bx,\dots,\tau^m(\bx),\by,\dots,\tau^m(\by))$$ and
$$d(tp(\bc,\dots,\tau^m(\bc),\bd,\dots,\tau^m(\bd)),tp(\bc',\dots,\tau^m(\bc'),\bd',\dots,\tau^m(\bd')))<\epsilon,$$
then $(\bc',\bd')\models
\psi'(\bx,\dots,\tau^m(\bx),\by,\dots,\tau^m(\by))$.

Let $b\in M$ and let $b'\in \F$ be such that
$d(tp(b',\dots,\tau^m(b')/A),tp(b,\dots,\tau^m(b)/A))<\epsilon$. Then
for all $\ba \subset A$, if $b\models \vphi(x,\ba)$ then
$b'\models \vphi'(x,\ba)$ and if $b'\models \vphi(x,\ba)$ then
$b\models \vphi'(x,\ba)$.
\end{proof}

\section{Ranks}

In this section we will follow \cite{Wal} and review the
definition and the main properties of entropy. Let $(X,\B,P)$ be a
probability space.

\begin{defi}
let $\A$ be a finite subalgebra of $\B$ with atoms
$\{A_1,\dots,A_k\}$. Let $\C$ be a sub-$\s$-algebra of $\B$. Then
the \emph{entropy of $\A$ given $\C$} is $$H(\A/\C)=-\int
\Sum_{i\leq k}\P(A_i|\C)\ln(\P(A_i|\C))dP$$
\end{defi}

We write $H(\A)$ for $H(\A/\{\emp,X\})$. If $\A$ and $\C$ are
$\s$-algebras, we denote by $\A \vee \C$ the $\s$-algebra
generated by $\A$ and $\C$.

\begin{fact}\label{properties0}
Let $\A$, $\C$ be finite subalgebras of $\B$ and let $\D$ be a
sub-$\s$-algebra of $\B$. Then:
\begin{enumerate}
\item $H(\A\vee\C/\D)=H(\A/\D)+H(\C/\A\vee \D)$.(Additivity)
\item $\A\subset \C \implies H(\A/\D)\leq H(\C/\D)$.
\item $H(\A/\D)\geq H(\A/C\vee \D)$,
\item If $\tau$ is a measure preserving automorphism, $H(\tau^{-1}\A/\tau^{-1}\D)=H(\A/\D)$
\item $H(\A/\C\vee \D)=H(\A/\D)$ iff $\A$ is independent from $\C$ over $\D$.
\end{enumerate}
\end{fact}

The first four properties are proved in Section 4.3 in \cite{Wal}.
A special case of property $(5)$, when $\D$ is the trivial algebra
$\{\emp,X\}$, is also proved there. Property $(5)$ was pointed out
by Ben-Yaacov to the authors and it can be proved using the ideas
presented in Section 4.3 in \cite{Wal}.

\begin{defi}
Let $\tau\colon X\to X$ be a measure preserving transformation
of the probability space $(X,\B,m)$. If $\A$ is a finite
subset of $\B$, then $$h(\tau,\A)=\lim_{n\rightarrow \infty}
\frac{1}{n}H(\bigvee_{i=0}^{n-1}\tau^{-i}\A)$$ is called the
\emph{entropy of $\tau$ with respect to $\A$}.

The value $h(\tau)=\sup\{h(\tau,\A):\A$ is a finite subset of $\B\}$ is
called the \emph{entropy} of $\tau$.
\end{defi}

\begin{fact}\label{properties}
Let $\A$ be a finite subalgebra of $\B$ and $\tau$ a measure
preserving automorphism of $(X,\B,m)$. Then:
\begin{enumerate}
\item $h(\tau,\A)\leq H(\A)$.
\item For $n>0$, $h(\tau,\A)=h(\tau,\bigvee_{i<n}\tau^{-i}\A)$.
\item $h(\tau,\A)=0$ iff $\A\subset \bigvee_{i=1}^{\infty}\tau^{-i}\A$.
\item $h(\tau,\A)=\lim_{n\rightarrow \infty}H(\A/\bigvee_{i=1}^{n}\tau^{-i}(\A))$.
\end{enumerate}
\end{fact}

The proofs can be found in Section 4.5 of \cite{Wal}.

\begin{defi}
Let $(X,\B,m)$ be an atomless probability space and let $\tau$ be
an aperiodic measure preserving automorphism of this space. Let
$M$ be the probability structure associated to $(X,\B,m)$. Let
$\ba=(a_1,\dots,a_n)\in M^n$ and let $D\subset M$. Let $A_1,\dots,A_n\in \B$
with events $a_1,\dots,a_n$ respectively and let $\A$ be the algebra
generated by $A_1,\dots,A_n$. Define
$H(\ba/D)$, the \emph{entropy of $\ba$ with respect to $D$}, to be
$H(\A/\bra D\ket)$. Define the \emph{entropy of $\tau$ with respect
to $\ba$} to be $h(\tau,\A)$ and denote it by $h(\tau,\ba)$.
Similarly, for $C=\{a_1,\dots,a_n\}$, let $H(C/D)=H(\A/\bra D\ket)$
and call it the \emph{entropy of $C$ with respect to $D$}. Finally
let the \emph{entropy of $\tau$ with respect to $C$} be
$h(\tau,\A)$ and denote it by $h(\tau,C)$.
\end{defi}

The properties listed in \ref{properties0}
and \ref{properties} still hold when measurable sets are replaced by
events.

\begin{defi}\label{transformally}
Let $(M,\tau)\models T_A$, let $\ba\in M^n$. We say that $\ba$ is
\emph{transformally independent} if $\{\tau^{i}(\ba):i\in
\mathbb{Z}\}$ is an independent sequence in $M$. We say that $\ba$
is \emph{transformally definable} if $\ba\in \dcl(\tau^{-i}(\ba):
i>0)$
\end{defi}

\begin{rema}
Let $(M,\tau)$ be a model of $T_A$, let $\ba\in M^n$. By the
additivity of entropy, $\ba$ is transformally independent if and
only if $h(\tau,\ba)=H(\ba)$. By property (3) in Fact
\ref{properties}, $\ba$ is transformally definable if and only if
$h(\tau,\ba)=0$. We call the algebra formed by the transformally
definable elements the \emph{Pinsker} algebra (compare with
\cite{Wal}).
\end{rema}

\begin{rema}
It is shown in Section 4.9 of \cite{Wal} that for every $r\in
\mathbb{R}^+$, there is a separable model $(M_r,\tau_r)\models
T_A$ such that $h(\tau_r)=r$. It follows that $h(\tau)=\infty$ whenever
 $(M,\tau)\models T_A$ is $\aleph_1$-saturated.
\end{rema}

\begin{obse}
We can use entropy to characterize the generic elements of the
groups $(M,\triangle)$, where $\triangle$ is the symmetric
difference. When we work in $M$ (just the probability structure),
it is shown in \cite{BY2} that the generics of the group are the
events of measure $1/2$. Note that if $a\in M$, then $H(a)\leq
\ln(2)$ and that the (maximal) value $\ln(2)$ is attained only
when $P(a)=1/2$. So the generics are the elements with maximal
entropy.

Similarly, if we work in the structure $(M,\tau)$, it is easy to
see that the generic elements of the group $(M,\tau,\triangle)$
are the transformally independent events of measure $1/2$. Let
$a\in M$. Then $h(a,\tau)\leq \ln(2)$ and equality holds iff $a$
is an event of measure $1/2$ which is transformally independent.
So the generics of $(M,\tau,\triangle)$ are the elements $a\in M$
such that the entropy of $\tau$ with respect to $a$ is maximal.
\end{obse}

\section{Orthogonality and omitting types}

We start with a definition from \cite{Wal}.

\begin{defi}
Let $(M,\tau)\models T_A$ and let $A\subset M$ be such that
$\tau(A)=A$. We say that $\tau$ has \emph{completely positive
entropy} on $A$ if for all finite $B\subset A$, $h(\tau,B)>0$.
\end{defi}

Let $M$ be a probability space structure and let $\tau$ be a
measure preserving aperiodic automorphism of the underlying
probability space. Let $A\subset M$ be countable. It is shown in
\cite{Pa} that $\tau$ has completely positive entropy on $A$ if
and only if there is $A'\subset A$ such that $\tau^{-1}(A')\subset
A'$, $\dcl(\cup_{n\geq 0} \tau^n(A'))\supset A$ and $\cap_{n\geq
1} \dcl(\tau^{-n}(A'))=\{\emp,X\}$.

It is well known (Theorem 4.37 in \cite{Wal}, Corollary
6.15 in \cite{Pa}) that if $\tau$ has completely positive entropy
on $A$, then $A$ is independent from the Pinsker $\s$-algebra of
$\tau$. A similar result is known about ACFA(\cite{Ma}), namely, 
that types of transformally independent
elements are orthogonal to types of transformally algebraic
elements. The combination of these two ideas suggests that
types of subsets of the Pinsker $\s$-algebra should be orthogonal to 
types of subsets where $\tau$ has completely positive entropy.
This result will be the first aim of this section. We need the
following notation from \cite{Pa}.

\begin{nota}
Let $(M,\tau)\models T_A$ and let $A\subset M$. Denote by $A^{-}$
the set $\dcl(\cup_{i\geq 1}\tau^{-i}(A))$.
\end{nota}

\begin{theo}
Let $(M,\tau)\models T_A$ be a $\kappa$-universal domain. Let
$A\subset M$ be a small set such that $\tau(A)= A$ and suppose
that $\tau$ has completely positive entropy on $A$. Let $B\subset
M$ be a small subset of the Pinsker $\s$-algebra. Let $\ba$ be an
enumeration of $A$ and $\bb$ an enumeration of $B$. Then
$tp_{\tau}(\ba)$ is orthogonal to $tp_{\tau}(\bb)$.
\end{theo}

\begin{proof}
We assume that $A$ and $B$ are countable. Then $A$ is
$\tau$-interdefinable with a set $A'$ such that
$\tau^{-1}(A')\subset A'$ and $\cap_{n\geq 1}
\dcl(\tau^{-n}(A'))=\{\emp,X\}$. Let $F\subset M$ be a small set
which is $\tau$-independent from $A'$ over $\emp$ and
$\tau$-independent from $B$ over $\emp$. Let
$F_\tau=\dcl_{\tau}(F)$ and $\ba'$ an enumeration of $A'$. We
need to show that
${\ba',\dots,\tau^n(\ba')}\ind_{F_\tau}{\bb,\dots.,\tau^m(\bb)}$ for
any $n,m<\omega$. Replacing $A'$ by $A'\cup\dots\cup \tau^n(A')$ and
$B$ by $B\cup\dots\cup \tau^m(B)$, it is enough to prove that
$\ba'\ind_{F_\tau}\bb$.

By the finite character of non-dividing, we only need to show that for
any finite subsets $C$ of $A'$ and $D$ of $B$, we have
${C}\ind_{F_\tau}{D}$.

By the additivity property of entropy, we can expand $H(C
\cup\dots\cup \tau^{-l}(C)\cup D \cup\dots\cup \tau^{-l}(D)/F_\tau)$ in two
different ways: $$H(C\cup \dots\cup \tau^{-l}(C)/F_\tau)+H(D\cup
\dots\cup \tau^{-l}(D)/F_\tau\cup C\cup \dots \cup \tau^{-l}(C))=$$
$$=H(D\cup \dots\cup \tau^{-l}(D)/F_\tau)+H(C\cup \dots\cup
\tau^{-l}(C)/F_\tau\cup D\cup \dots \cup \tau^{-l}(D))$$ Dividing by
$l+1$, taking limits as $l$ goes to infinity and using Fact
\ref{properties}, we get $$H(C/F_\tau\cup C^-)+H(D/F_\tau\cup
D^-\cup C^-\cup C)=H(D/F_\tau\cup D^-)+H(C/F_\tau\cup C^-\cup D
\cup D^-)$$

Since $D$ is a subset of the Pinsker $\s$-algebra, this implies
$H(C/F_\tau \cup C^-)=H(C/F_\tau \cup C^- \cup D \cup D^-)$. Thus
${C}\ind_{F_\tau \cup C^-}{D\cup D^{-}}$. The same argument holds
if we replace $C$ by $\tau^{-i}(C)$, so we get
${\tau^{-i}(C)}\ind_{F_\tau \cup (\tau^{-i}(C))^-}{D\cup D^{-}}$
for all $i\geq 0$. By transitivity of independence, this implies
${C\dots\tau^{-i}(C)}\ind_{F_\tau \cup (\tau^{-i}(C))^-}{D\cup
D^{-}}$ for any $i\geq 0$. Since $\cap_{i\geq 1}
\dcl(\tau^{-i}(A'))=\{\emp,X\}$, we get ${C}\ind_{F_\tau}{D\cup
D^{-}}$ as we wanted.
\end{proof}

Since there are $2^{\aleph_0}$ many non-isomorphic Bernoulli
shifts (see Section 4.9 in \cite{Wal}), all of which induce aperiodic
maps on separable complete probability structures, there are many 
non principal types over $\emp$. In the rest of this section
we prove the stronger result that only algebraic types are principal. 
We start with some definitions.

\begin{defi}
Let $(M,\tau)\models T_A$ and let $B\subset M$ be such that
$\tau(B)=B$. If $A\subset M$ is a finite set, then
$$h(\tau,A/B)=\lim_{n\rightarrow \infty}
\frac{1}{n}H(\bigcup_{i=0}^{n-1}\tau^{-i}A/B)$$ is called the
\emph{entropy of $\tau$ with respect to $A$ over $B$}.

We also define $h(\tau /B)=\sup\{h(\tau,A/B):A$ is a finite subset of $M\}$. This is called the \emph{entropy of $\tau$ with respect to
$B$}.
\end{defi}

\begin{obse}\label{finite}
Let $(M,\tau)\models T_A$ and let $B\subset M$ be such that
$\tau(B)=B$.

\begin{enumerate}
\item Let $A$ be a finite subalgebra of $M$. Then
$h(\tau,A/B)=H(A/B\cup A^-)$.
\item Let $C$ be a subalgebra of $M$ such that $\dcl(C)=M$.
Then $h(\tau/B)=\sup\{h(\tau,A/B):A$ is a finite subset of $C\}$.
\end{enumerate}
\end{obse}

\begin{proof}
When $B=\emp$, (1) is proved in Section 4.5 in \cite{Wal} and (2)
is proved in Section 4.6 in \cite{Wal}. The proofs in \cite{Wal}
easily generalize to give what is stated here.
\end{proof}

\begin{prop}\label{omit0}
Let $(M,\tau)\models T_A$ be a $\kappa$-universal domain and let
$B\subset M$ be countable such that $\tau(B)=B$. Then there is
$(M_1,\tau_1)\subset(M,\tau)$ such that $(M_1,\tau_1)\models T_A$,
$B\subset M_1$ and for any finite algebra $D\subset M_1$,
$h(\tau,D/B)=0$.
\end{prop}

\begin{proof}
Let $(M_0,\tau_0)$ be the model of $T_A$ induced by an
irrational rotation on the unit circle. Then $h(\tau,D)=0$ for any finite
subalgebra $D\subset M_0$ (see Section 4.7 in
\cite{Wal}). Since $(M,\tau)$ is $\kappa$-saturated, we may 
assume that  $(M_0,\tau_0)$ is a substructure of $(M,\tau)$ which is 
independent from $B$. Let $M_1=\dcl(M_0,B)$ and let $\tau_1$
 be the restriction of $\tau$ to $M_1$.

Let $B'=\dcl(B)$. Then $(B',\tau \upharpoonright_{B'})$ is
separable and can be seen as a probability structure with an 
automorphism. The
structure $(M_1,\tau_1)$ is isomorphic to the product of the
structures $(M_0,\tau_0)$ and $(B', \tau\upharpoonright_{B'})$.
Since $\tau_0$ is an aperiodic map, we get $(M_1,\tau_1)\models
T_A$.

It remains to show that for any finite algebra $D\subset M_1$,
$h(\tau,D/B)=0$. By \ref{finite}, it is enough to prove
that for any finite algebras $B_0\subset B$ and $A_0\subset M_0'$,
$h(\tau,B_0 \cup A_0)/B)=0$. Now, since the entropy of $\tau_0$ is
zero, the following inequalities hold: $h(\tau,B_0 \cup
A_0/B)=H(B_0\vee A_0/B\cup A_0^-)\leq H(A_0/A_0^-)=0$.
\end{proof}

\begin{prop}\label{omit1}
Let $(M,\tau)\models T_A$ be a $\kappa$-universal domain and let
$B\subset M$ be countable such that $\tau(B)=B$. Then there is
$(M_1,\tau_1)\subset (M,\tau)$ such that $(M_1,\tau_1)\models
T_A$, $B\subset M_1$ and for any finite algebra $D\subset M$,
$h(\tau,D/B)=0$ iff $D\subset \dcl(B)$.
\end{prop}

\begin{proof}
Let $(M_0,\tau_0)\models T_A$ be the structure induced by a
Bernoulli shift generated by a partition into two elements
$\{a_1,a_2\}$ of probability $1/2$ each. Then for any finite set
$D\subset M_0$, $h(\tau,D)=0$ iff $D\subset \{\emp,X\}$ (see
Section 4.9 in \cite{Wal} or \cite{Sh}).

Since $(M,\tau)$ is $\kappa$-saturated, we may assume
$(M_0,\tau_0)$ is a substructure of $(M,\tau)$ which is
independent from $B$. Let $B'=\dcl(B)$. Then $(B',\tau
\upharpoonright_{B'})$ is separable. Let $M_1=\dcl(M_0,B)$ and let
$\tau_1$ be the restriction of $\tau$ to $M_1$. Then
$(M_1,\tau_1)$ is isomorphic to the product of $(M_0,\tau_0)$ and
$(B',\tau \upharpoonright_{B'})$. Since $\tau_0$ is an aperiodic
map, then $(M_1,\tau_1)\models T_A$.

Let $D\subset M_1$ be a finite subalgebra and assume that
$h(\tau,D/B)=0$. We want to show that $D\subset \dcl(B)$.

Let $A=\{a_1,a_2\}$. Let $A^-=\dcl(\cup_{i\geq 1}\tau^{-i}(A))$
and $D^-=\dcl(\cup_{i\geq 1}\tau^{-i}(D))$. Then by the additivity
property of entropy, we can show that $H(A/ A^- \cup B\cup D \cup
D^-)+H(D/B \cup D^-)= H(D/ D^- \cup B\cup A\cup A^-)+H(A/B \cup
A^-)$. Since $h(\tau,D/B)=0$, we get $H(A/ A^- \cup B \cup D\cup
D^-)=H(A/B \cup A^-)$. This proves ${A}\ind_{B\cup
A^{-}}{D}$. Exchanging $A$ for $\tau^{m}(A)$, we get
${\tau^m(A)}\ind_{B\cup (\tau^m(A))^-}{D}$ for any $m\in
\mathbb{Z}$. By transitivity of independence we obtain
${\tau^m(A),\dots,\tau^{-m}(A)}\ind_{B\cup (\tau^{-m}(A))^-}{D}$ for
any $m\in \mathbb{N}$. Since $A$ is transformally independent and
$\tau$-independent from $B$, we obtain
${\tau^m(A),\dots,\tau^{-m}(A)}\ind_{B}{D}$ for any $m\in
\mathbb{N}$. By the finite character of independence, this implies
${M_0}\ind_{B}{D}$. Since $D\subset \dcl(M_0,B)$, we must also
have $D\subset \dcl(B)$.
\end{proof}

\begin{prop}
Let $(M,\tau)\models T_A$ be a $\kappa$-universal domain and let
$B\subset M$ be countable such that $\tau(B)=B$. Let $\ba\in M^n$.
Then $tp(\ba/B)$ is principal iff $\ba\in \dcl(B)^n$.
\end{prop}

\begin{proof}
Let $B$ and $\ba$ be as above and assume that $\ba \not \in \dcl(B)^n$. 
If $h(\tau,\ba/B)>0$, then we can
omit $tp(\ba/B)$ by Proposition \ref{omit0}. If $h(\tau,\ba/B)=0$,
then we can omit $tp(\ba/B)$ by Proposition \ref{omit1}.
\end{proof}

\section{General automorphisms}

Let $(X,\B,m)$ be a probability space, let $\tau$ be an
automorphism of this space and let $M$ be the probability
structure associated to $(X,\B,m)$. The aim of this section is to
discuss $Th_{\A}(M,\tau)$.

Let $Y$ be the union of the atoms from $(X,\B,m)$ and let $\B_Y$ be
the $\s$-algebra induced by $\B$ on $Y$. Then $\tau$ is an
automorphism of $Y$ that permutes each set of atoms having the same
(positive) measure. First we characterize $((Y,\B_Y,m),\tau)$. Assume
that $((Y',\B',m'),\tau')$ is another structure such that
$(\events(Y',\B',m'),\tau')\equiv_\A (\events(Y,\B_Y,m),\tau)$. Then for every
real number $r\in (0,1]$, the number of atoms in $\B_Y$ with
measure $r$ (denoted by $A_r$) agrees with the number of atoms of
$\B'$ with measure $r$ (denoted by $A_r'$) and the action of
$\tau$ on $A_r$ is isomorphic to the action of $\tau'$ on $A_r'$.
So for every $r\in (0,1]$, $(A_r,\tau)\cong (A_r',\tau')$. This
implies $(\events(Y',\B',m'),\tau')\cong (\events(Y,\B_Y,m),\tau)$. 
Thus, to characterize $(\events(Y,\B_Y,m),\tau)$, we only need to 
describe the permutation that $\tau$ induces on the set of atoms of measure $r$ for each $r\in (0,1]$. The structure $(\events(Y,\B_Y,m),\tau)$ behaves like
a finite structure in first order theories.

Let $Z$ be the atomless part of $X$. We can decompose $Z$
into a disjoint union $\cup_{i\in \mathbb{N}}Z_i$, where $Z_0$ is
the set of aperiodic elements of $Z$ and for $i>0$, $Z_i=\{x\in X:
\tau^i(x)=x\}\setminus (Z_1\cup\dots\cup Z_{i-1})$. The automorphism
$\tau$ acts on each of the sets $Z_i$. Let $\B_{Z_i}$ be the
$\s$-algebra induced by $\B$ on $Z_i$. Let $M_i$ be the
probability structure associated to $(Z_i,\B_{Z_i},m)$. To study
the atomless part of $Th_{\A}(M,\tau)$, it suffices to understand
$Th_{\A}(M_i,\tau)$ for $i\in \mathbb{N}$. The behavior of the
aperiodic part $Th_{\A}(M_0,\tau)$, is
described by $T_A$ after rescaling $m(Z_0)$ to be $1$.

\begin{lema}\label{partition}
Let $([0,1],\B,m)$ be the standard Lebesgue space, let $n\in
\mathbb{N}$ and let $\tau$ be an automorphism such that
$\tau^{n+1}=id$ and $m(\{x:\tau^j(x)=x\})=0$ for all $j<n+1$. 
Then there is a set $A\in \B$ such that $(A,
\dots,\tau^{n}(A))$ forms a partition of $[0,1]$ up to measure zero.
\end{lema}

\begin{proof}
See Lemma 1 in \cite[pp. 70]{Ha}.
\end{proof}

We can now show that $Th_{\A}(M_i,\tau)$ is separably categorical 
for each $i\geq 1$.

\begin{prop}\label{sepcat}
Let $([0,1],\B,m)$ be the standard Lebesgue space, let $i\in
\mathbb{N}^+$ and let $\tau,\eta$ be an automorphisms such that
$\tau^i=id$, $\eta^i=id$ and $m(\{x:\tau^j(x)=x\})=0$,
$m(\{x:\eta^j(x)=x\})=0$ for all $j<i$. Let $N$ be the probability
structure associated to $([0,1],\B,m)$. Then $(N,\tau)\cong
(N,\eta)$.
\end{prop}

\begin{proof}
By the previous lemma, there are $a\in N$ and $b\in N$ such that
$(a,\tau(a),\dots,\tau^{i-1}(a))$ forms a partition of $N$ and
$(b,\eta(b),\dots,\eta^{i-1}(a))$ forms a partition of $N$. Let
$A\in \B$ with event $a$ and let $B\in \B$ with event $b$. There
is a measure preserving automorphism $\colon:A\to B$. We
can extend $\gamma$ by defining for $x=\tau^i(y)\in \tau^i(A)$, 
$\gamma(x)=\eta^i(\gamma(y))$. Then $\gamma \tau=\eta \gamma$ 
on a set of measure one. This proves that $(N,\tau) \cong (N,\eta)$.
\end{proof}

\begin{lema} Let $([0,1],\B,m)$ be the standard Lebesgue space, 
let $n\in \mathbb{N}^+$ and let $\tau$ be an automorphism such 
that $\tau^{n+1}=id$ and $m(\{x:\tau^j(x)=x\})=0$ for all $j<n+1$. 
Let $B_1,\dots,B_m\in \B$. Then there is a set $A\in \B$ of measure 
$1/(n+1)$ such that $(A,\dots,\tau^{n}(A))$ forms a partition of 
$[0,1]$ (up to measure zero) which is independent from $\{B_1,\dots,
B_m\}$.
\end{lema}

\begin{proof}
Replacing $B_1,\dots,B_m$ by the atoms in the algebra generated 
by $\{\tau^i(B_j): 0\leq i \leq n, 1\leq j\leq m\}$, we may assume 
that $\tau$ acts on the set $\{B_1,\dots,B_m\}$. If $\tau(B_i)=B_i$, by Lemma 
\ref{partition}, we can find $A_i\in \B$ such that $A_i\subset B_i$, 
and $(A_i,\dots,\eta^{n}(A_i))$ forms a partition of $B_i$. If $\tau$ 
acts transitively on $B_{i_1},\dots,\B_{i_k}$, just repeat the argument 
from \cite[pp 70]{Ha} starting with $B_{i_1}$ instead of $E_1$. From the 
sets $\{A_i:1\leq i\leq m\}$ we can construct $A$.
\end{proof}

Now we prove that $Th_{\A}(M_i,\tau)$ has quantifier elimination for each
$i\geq 1$.

\begin{lema}\label{qelim}
Let $([0,1],\B,m)$ be the standard atomless Lebesgue space and let
$M$ be the probability structure associate to it. Let $\eta_1$,
$\eta_2$ be cycles of period $n+1$ on $M$. Let $\bb,\bd\in M^m$ be 
such that $qftp(\bb,\dots,\eta_1^n(\bb))=qftp(\bd,\dots,\eta_2^n(\bd))$.
Then there is an automorphism $\gamma$ of $M$ such that $\gamma
\eta_1=\eta_2 \gamma$ and $\gamma(\bb)=(\bd)$.
\end{lema}

\begin{proof}
By the previous lemma we can find $a\in M$ such that $(a,\eta_1(a),\dots,
\eta_1^n(a))$ is a partition of $M$  
which is independent from $\{\bb,\dots,\eta_1^n(\bb)\}$.

Similarly there is $c\in M$ such that $(c,\eta_2(c),\dots,
\eta_2^n(c))$ is a partition of $M$ which is independent 
from $\{\bd,\dots,\eta_2^n(\bd)\}$.

In particular we obtain that $qftp(\bb,\dots,\eta_1^n(\bb),a,\dots,
\eta_1^n(a))=qftp(\bd,\dots,\eta_2^n(\bd),c,\dots,\eta_2^n(c))$.
Let $\bb'=(b_1',\dots,b_{k}')$ be the atoms of the algebra generated 
by  $\{\bb,\dots,\eta_1^n(\bb)\}$ and let $\bd'=(d_1',\dots,d_{k}')$ be 
the atoms of the algebra generated by $\{\bd,\dots,\eta_2^n(\bd)\}$. We may
choose enumerations such that $qftp(\bb',a,\dots,\eta_1^n(a))=
qftp(\bd',c,\dots,\eta_1^n(c))$. It suffices to prove the lemma 
for $\bb'$ and $\bd'$.

Let $A,B_1',\dots,B_k',C,D_1',\dots,D_k'$ be sets in $\B$ giving rise to 
the collection of events $a,b_1',\dots,b_k',c_1,\dots,d_k'$. Then $m(A\cap B_i')=m(C\cap D_i')$ for $i\leq k$. Let
$\gamma \colon A\to C$ be an isomorphism sending $A\cap B_i$ to
$C\cap D_i$ (up to measure zero) for $i\leq k$. Let $x\in
\eta_1^l(A)$. Then there is $y\in A$ such that $\eta_1^l(y)=x$.
Extend $\gamma$ by defining $\gamma(x)=\eta_2^l(\gamma(y))$. Note
that $\eta_1(b_i')=d_i'$ for $i\leq k$ and $\gamma
\eta_1=\eta_2 \gamma$.
\end{proof}

By Lemma \ref{qelim}, $Th_{\A}(M_i,\tau)$ has quantifier
elimination and thus it is model complete.
Note that $Th_{\A}(M,\tau)$ need not have quantifier elimination.
We first need to split $(M,\tau)$ into its atomic part and the $(M_i,\tau)$
 in order to get quantifier elimination.

Since $T$ is stable, so are the theories $Th_{\A}(M_i,\tau)$ for
$i>0$ and $Th_{\A}(M,\tau)$.

Let $(M,\tau)\preceq_{\A} (N,\tau)$ and let $\{z_i:i \in
\mathbb{N}\}$ be the events corresponding to the measurable sets
$\{Z_i :i \in \mathbb{N}\}$ described above. Let $i\in
\mathbb{N}^+$ and assume that $m(z_i)>0$. Let $a\in N$ be such
that $a\subset z_i$ and $m(a)>0$. Then there is $b\subset a$ such
that $m(b)>0$, $\tau^i(b)=b$ and $\tau^j(b)\neq b$ for all $j<i$.
This shows the decomposition $\{z_i:i \in \mathbb{N}\}$ is
preserved in elementary superstructures.

Let $(X',\B',m')$ be a probability space, let $\tau'$ be a measure
preserving automorphism of this space and let $M'$ be the
probability structure associated to the probability space. Follow
the notation above and denote by $Y'$ the union of the atoms
from $M'$ and $z_i'$ the event associated to the set formed by the 
elements with period $i$. Then $(M',\tau')\equiv_{\A} (M,\tau)$ iff
$((Y',\B_Y',m'),\tau')\cong ((Y,\B_Y,m),\tau)$ and
$m(z_i')=m(z_i)$ for $i\in \mathbb{N}$. We conclude that the
theory $Th_{\A}(M,\tau)$ can be described in terms of
$Th_{\A}(M_Y,\tau)$ and the sequence $(m(M_i):i<\omega)$.

In conclusion, we give some open questions related to the
subject of this paper.

As background to the first questions, note that for
probability structures (\textit{i.e.}, without
automorphism), Lemma \ref{distance} gives a very nice,
explicit formula for the $d$ distance between types over a
set of parameters $C$.  This is a basis for a full analysis
of stability of these structures.  For example, see Remark
\ref{superstable}.

\begin{ques}
Can we characterize the $d$-metric in type spaces of $T_A$?
That is, can we find a description of $d$ in the spirit of
Lemma \ref{distance}?
\end{ques}

\begin{ques}
What is the density character of the space of $\tau$-types
over a given set of parameters $C$, with respect to the
$d$-metric?  (See \cite{Io1,Io2} for a proof that $T_A$ is
stable with respect to the $d$-metric, since it is stable
with respect to the minimal uniform structure on types.
This gives a partial answer to this question.)
\end{ques}

\begin{ques} 
Is $T_A$ superstable with respect to the
$d$-metric? (See Remark \ref{superstable}.)
\end{ques}

In \cite{BBH} the model theory of the Banach lattices $L_p(\mu)$
is studied, for each $p \in [1,\infty)$.  Each of these theories
interprets the theory $T$ of probability structures (for any
positive $L_p$-function $f$ of norm 1, consider the set of components
of $f$ equipped with the $p^{th}$ power of the norm as probability
measure).  So, in a certain sense the results in \cite{BBH} extend
those in \cite{BY2} and in section 3 of this paper.  Furthermore,
automorphisms of $L_p(\mu)$-spaces are well understood from the 
functional analysis point of view.  

\begin{ques}
Is there a model companion for the positive theory of the Banach
lattices $L_p(\mu)$ expanded by an automorphism?
\end{ques}

\end{document}